\documentclass[12pt]{amsart}

\usepackage{amssymb}
\usepackage{amsmath}
\usepackage{amsthm}

\input xy
\xyoption{all}

\title{Universal Operator Algebras of Directed Graphs}

\author{Benton L. Duncan}

\address{Department of Mathematics\\
810 Oldfather Hall\\
University of Nebraska, Lincoln\\
Lincoln NE  68588-0323}

\email{bduncan@math.unl.edu}



\thanks{Part of this work was supported by a Department of Education GAANN fellowship}

\begin{document}

\theoremstyle{plain}
\newtheorem{theorem}{Theorem}[section]
\newtheorem{lemma}{Lemma}[section]
\newtheorem{proposition}{Proposition}[section]
\newtheorem{corollary}{Corollary}[section]
\newtheorem{conjecture}{Conjecture}[section]

\theoremstyle{definition}
\newtheorem{definition}{Definition}[section]
\newtheorem{construction}{Construction}[section]
\newtheorem{example}{Example}[section]

\theoremstyle{remark}
\newtheorem{remark}{Remark}[section]
\newtheorem{question}{Question}
\newtheorem*{acknowledgement}{Acknowledgements}

\begin{abstract} Given a directed graph, there exists a universal
operator algebra and universal $C^*$-algebra associated to the
directed graph.  For finite graphs this algebra decomposes as the
universal free product of some building block operator algebras.
For countable directed graphs, the universal operator algebras
arise as direct limits of operator algebras of finite subgraphs.
Finally, a method for computing the K-groups for universal
operator algebras of directed graphs is given.
\end{abstract}

\maketitle

In \cite{Muhly:1997} Muhly associates a non-selfadjoint operator
algebra to a directed graph (or quiver).  Henceforth we refer to
these algebras as Toeplitz quiver algebras.  Kribs and Power
\cite{Kribs-Power:2003a} showed that the graph was a complete
unitary invariant for these algebras.  Recent work on these
Toeplitz quiver algebras by Katsoulis and Kribs,
\cite{Kat-Kribs:2003} and by Solel \cite{Sol:2003}, has
demonstrated that the graph is a complete isomorphism invariant
for these algebras. In addition Kribs and Power,
\cite{Kribs-Power:2003a} and \cite{Kribs-Power:2003b} study the
structure of these algebras, including a free product result for
certain amalgamations of graphs.

In another direction \cite{Bl-Paul:1991}, and \cite{Good-Men:1990}
have initiated a study of universal operator algebras (both
nonself-adjoint and self-adjoint) associated to combinatorial
objects (e.g. groups, monoids, and semigroups). We have continued
this study by introducing the universal operator algebra, and the
universal $C^*$-algebra associated to a directed graph.  In what
follows we look at the construction of these objects, and using
\cite{Good-Men:1990}, \cite{Kribs-Power:2003a}, and
\cite{Spielberg:2002} as models, we find nice decompositions. As a
consequence of the decompositions we calculate the K-Theory of
these algebras.

As a result of theorem \ref{freeproducts} we are able to write any
finite graph as a free product of copies of three ``building block
graphs''. Using universal properties, we then show that the
universal operator algebra will be a free product of three
building block algebras. For the universal $C^*$-algebra of a
directed graph we use a construction of Blecher \cite{Bl:1999} to
also decompose the universal $C^*$-algebra of a directed graph as
free products.

For infinite graphs we are able, as in \cite{Spielberg:2002}, to
write the universal operator algebra as a direct limit of the
universal operator algebras corresponding to finite subgraphs
directed by inclusion.  We then prove that the maximal
$C^*$-algebra of an operator algebra preserves direct limits, and
hence the same result for universal $C^*$-algebras of directed
graphs follows.

These decompositions, with a result of Cuntz \cite{Cu:1982}, allow
us to compute the K-groups for the universal operator algebras of
directed graphs.  We are then able to extend our calculations to
the K-groups of the universal $C^*$-algebra of a directed graph.
Our calculations of the K-groups indicates that they are a poor
invariant in studying these non-selfadjoint operator algebras.  In
fact, for a directed graph the K-groups depend only on the number
of vertices.

Looking at our results on a categorical level we have associated
two functors on the category of directed graphs (with directed
graph homomorphisms): the first functor is into the category of
operator algebras (with completely contractive homomorphisms); the
second functor is into the category of $C^*$-algebras (with $C^*$
homomorphisms).  We show that the two functors behave well with
respect to free products and direct limits.  We notice that the
functors presented here differ from those in
\cite{Spielberg:2002}, since the morphisms he considers are
directed graph inclusions.

Before proceeding, we would like to emphasize a difference between
the operator algebras in the present paper and those defined in
\cite{Muhly:1997}. When we construct the universal operator
algebra of a directed graph we consider representations which send
vertices to projections.  \emph{We do not assume that the
projections are orthogonal}, as was implicit in \cite{Muhly:1997}.
This difference provides examples which differ significantly from
the Toeplitz quiver algebras.  It is also important to distinguish
here between the universal free products used herein, and the
``spatial'' free products in \cite{Kribs-Power:2003a}. Kribs and
Power construct a specific representation of the free product when
proving a decomposition result for directed graphs. Our
decomposition, on the other hand, takes into account all
representations of the two graphs, as in \cite{Bl-Paul:1991}. This
yields a significant distinction, even in the case of the directed
graph given by two vertices with no edges.

\section{Preliminaries}\label{prelims}

In \cite{Bl-Paul:1991}, Blecher and Paulsen introduced the maximal
$C^*$-algebra of an operator algebra.  For $A$ an operator algebra
recall, as in \cite{Bl:1999}, the construction of the maximal
$C^*$-algebra, $C^*_m(A)$.

\begin{theorem}[Blecher \cite{Bl:1999} Theorem 2.1]\label{cm} For a unital
operator algebra $A$ there exists a unital $C^*$-algebra,
$C^*_m(A)$ and a completely isometric inclusion $ \iota: A
\rightarrow C^*_m(A)$ such that: \begin{enumerate} \item $
\iota(A) $ generates $C^*_m(A)$ as a $C^*$-algebra.
\item\label{cmup} (Universal Property) for any unital completely
contractive homomorphism $\varphi: A \rightarrow C$ into a
$C^*$-algebra $C$ there is a unique $*$-homomorphism $
\tilde{\varphi}: C^*_m(A) \rightarrow C$ making the following
diagram commute
\[
\xymatrix{ C^*_m(A) \ar@{-->}[dr]^{\tilde{\varphi}} & \\ A
\ar[u]^{\iota} \ar[r]^{\varphi} & C }
\]
\end{enumerate}
\end{theorem}
\begin{remark} If $A$ is nonunital we denote the unitization
of $A$ by $A^+$.  By the universal property of the unitization,
\cite[Theorem 3.2]{Mey:2001}, $C^*_m(A)$ is the $C^*$ subalgebra
of $C^*_m(A^+)$ generated by $A$.\end{remark} Notice that by the
universal property (\ref{cmup}) of Theorem \ref{cm}, $C^*_m(A)$ is
unique up to $*$-isomorphism.

Before the construction, we recall some facts about operator
algebras that can be found in \cite{Paul:2002}. For an operator
algebra $A$, the associated operator algebra $A^*$ is unique up to
completely isometric isomorphism, independent of the particular
representation of $A$. Further notice  that for any completely
contractive homomorphism $ \varphi: A \rightarrow C$ into a
$C^*$-algebra there is, by Arveson's extension formula
\cite[Theorem 1.29]{Arv:1969}, a unique completely contractive map
$ \overline{\varphi}:A + A^* \rightarrow C$, which extends
$\varphi$ and such that $\overline{\varphi}(a^*) = \varphi(a)^*$
for $a \in A$. We will denote $ \overline{\varphi}|_{A^*}$ by
$\varphi^*$.

We now go through the construction of the maximal $C^*$-algebra as
it is important in what follows.
\begin{construction}\cite[Theorem 2.1]{Bl:1999} First form the algebraic free product
$ A *_a A^*$ of $A$ and $A^*$ amalgamated over the diagonal $A
\cap A^*$. This is a $*$-algebra, where $(a_1*a_2*\cdots*a_n)^* =
(a_n^* * a_{n-1}^* * \cdots * a_1^*)$ and the involution is
extended using linearity. By the universal property for $A*_aA^*$,
any completely contractive representation $ \varphi :A \rightarrow
B(H_{\varphi})$ extends to a $*$-representation, $\varphi *
\varphi^*$ of $A*_aA^*$.  Now define \begin{align*} \| x \| = \sup
\{ \| \varphi * \varphi^* (x) \|_{B(H_{\varphi})} &: \varphi
\mbox{ is a completely} \\ & \mbox{contractive representation of }
A \} \end{align*} for all $ x \in A*_aA^*$. Completing $A*_aA^*$
with respect to this norm will give a $C^*$-algebra satisfying the
universal properties given in Theorem \ref{cm} (\ref{cmgup}).
\end{construction}
\begin{remark}The construction shows that \[ P = \mbox{linspan} \cup_{n
\in \mathbb{N}}\{ a_1*a_2*\cdots *a_n : a_i \in A \cup A^*\}\]
forms a dense set in $C^*_m(A)$.\end{remark}

\section{The Universal Operator Algebra of a Directed Graph}

Now we let $Q$ be a directed graph with vertex set $V(Q)$ and edge
set $E(Q)$.  To each edge $e$ we denote by $s(e)$ the source
vertex for $E$ and by $r(e)$ the range vertex for $e$.  For $n
\geq 1$ we let $E(Q)^{n}$ be the set of words over the edges of
length $n$, and we let $E(Q)^0 $ be the set of words over the
vertex set. \emph{We allow arbitrary words, with no restriction to
admissible words, as in \cite{Muhly:1997}.} We define
\[ w(Q) = \bigcup_{n=0}^{\infty} E(Q)^n. \] We denote, by $|w|$,
the length of the word $ w \in w(Q)$ (i.e $|w | = n \mbox{ if } w
\in E(Q)^n$).

\begin{definition} For an operator algebra $A$ we say that $\pi: Q \rightarrow A$
is a contractive representation of $Q$ if \begin{enumerate} \item
$\pi(v)$ is a projection for all $v \in V(Q)$. \item $\|\pi(e)\|
\leq 1$ for all $ e \in E(Q)$. \item $\pi(v)\pi(e) = \pi(e)$ if
$r(e) = v$ and $\pi(e) \pi(v) = \pi(e) $ if $s(e) = v$.
\end{enumerate}\end{definition}
\begin{proposition}\label{oag} For a directed graph $Q$ there exists
a unique operator algebra $OA(Q)$ and a contractive representation
$\iota: Q \rightarrow OA(Q)$ such that:
\begin{enumerate} \item $ \iota(Q) $ generates $OA(Q)$ as an
operator algebra. \item\label{oagup} (Universal Property) When $A$
is an operator algebra and $ \varphi: Q \rightarrow A$ is a
contractive representation there is a unique completely
contractive homomorphism $ \tilde{\varphi}: OA(Q) \rightarrow A$
making the following diagram commute
\[ \xymatrix{ OA(Q)
\ar@{-->}[dr]^{\tilde{\varphi}} & \\ Q \ar[u]^{\iota}
\ar[r]^{\varphi} & A } \]
\end{enumerate}
\end{proposition}
\begin{remark} Notice that $OA(Q)$ is unique up to completely
isometric isomorphism. \end{remark}
\begin{proof}  Let $\mathbb{C}w(Q)$ be the set of complex valued
functions on $w(Q)$ with finite support. Define a multiplication
on $\mathbb{C}w(Q)$ by \[ f*g(w) = \sum_{w_1w_2 = w} f(w_1)g(w_2).
\] Then $\mathbb{C}w(Q)$ is an algebra with pointwise
addition and multiplication as defined.

Any contractive representation $ \pi: Q \rightarrow A$ extends
uniquely to a representation $\tilde{\pi}: \mathbb{C}w(Q)
\rightarrow A$, by sending \[ f(w) \mapsto f(w)\pi(e_1) \pi(e_2)
\cdots \pi(e_n)\] where $w = e_1e_2 \cdots e_n$.

We define
\[ \| f \|_{OA} = \sup \{ \| \tilde{\pi}(f) \|: \pi \mbox{ is a
contractive representation of } Q\}.\] Notice that $\|
\tilde{\pi}(f) \|_A \leq \| f \|_1< \infty $ for a contractive
representation $\pi$, and hence $\| f \|_{OA} < \infty$.

To verify that this is a norm we need to show that for $f \neq 0
\in \mathbb{C}w(Q)$ there exists a representation $\pi: Q
\rightarrow A$ such that $\tilde{\pi}(f) \neq 0$.  Let $w \in
w(Q)$ be chosen so that \[ |w| = \min\{ |x| : x \in \mbox{supp}(f)
\}. \]  If $ |w| = 0$, then $w = v_1 v_2 \cdots v_k$ is a product
of $k$ vertices for some $k$. Let $\mathbb{C}_1* \mathbb{C}_2 *
\cdots * \mathbb{C}_k$ be the nonamalgamated free product of $k$
copies of $ \mathbb{C}$ and define $\pi: Q \rightarrow
\mathbb{C}_1* \mathbb{C}_2 * \cdots \mathbb{C}_k$ by
\[ \pi(v_i) =
\begin{cases}  1_{i} & \mbox{ where }1_i \mbox{ is the unit
in the }i\mbox{th copy of }\mathbb{C} \\ 0 & \mbox{ else.}
\end{cases} \]  Then $ \tilde{\pi}(f) \neq 0$.

If $ | w | > 0 $, let $ n = |w|$ with $ w = e_1e_2 \cdots e_n$ and
consider the left regular representation of $A_n$, the
noncommmutative analytic Toeplitz algebra, (see
\cite{Dav-Pitts:1998} or \cite{Pop:1998}), with generators $T_1,
T_2, \cdots, T_n$. Now define $\pi: Q
\rightarrow A_n$ via \[ \pi(x) = \begin{cases} 1 & x \in V(Q) \\
T_i & x = e_i
\\ 0 & \mbox{ else.} \end{cases} \] Notice first that $ \pi$ is
a contractive representation of $Q$. Secondly notice that, by the
definition of the left regular representation, see the discussion
in section 6 of \cite{Muhly:1997}, we have a faithful
representation of $A_n$, and hence $\tilde{\pi}(f) \neq 0$. It
follows that $\| \cdot \|_{OA} $ is indeed a norm.

Completing $\mathbb{C}w(Q)$ with respect to this norm yields an
operator algebra satisfying the requisite universal property.
\end{proof}
\begin{remark} Notice that we are not claiming that $A_n$ is
universal in the sense of Theorem \ref{oag}.  In fact, it will
fall out of our later analysis that this is not the case. $A_n$ is
just used as a tool to verify that $\| \cdot \|_{OA}$ is a
norm.\end{remark}

\begin{example}\label{dot} Let $V_1$ be the finite graph given by
\[
\xymatrix{ {\bullet}_{v_0}}
\]
Here the universal operator algebra is $\mathbb{C}$.
\end{example}
\begin{example}\label{doubledot} Let $V_2$ be the finite graph given by
\[
\xymatrix{ {\bullet}_{v_0} & {\bullet}_{v_1} }
\]
Somewhat surprisingly $OA(V_2)$ is not isomorphic to $\mathbb{C}
\oplus \mathbb{C}$.  In fact, this algebra turns out to be the
unamalgamated universal free product $ \mathbb{C} * \mathbb{C}$.
\end{example}
\begin{example}\label{loop} We next look at the graph $B_1$ given by:
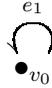
\[
\xymatrix{ {\bullet}_{v_0} \ar@(ur,ul)[]_{e_1} }
\]
We claim that $OA(B_1)$ is completely isometrically isomorphic to
the disk algebra, $A(\mathbb{D})$.  First define the contractive
representation $\iota: B_1 \rightarrow A(\mathbb{D})$ by sending
$v_0$ to the identity and $e_1$ to the coordinate function $f(z) =
z$. Notice also that $\iota(B_1) $ generates $A(\mathbb{D})$.

Now if $ \pi: B_1 \rightarrow A$ is a contractive representation
of $B_1$, then define the contractive homomorphism $\tilde{\pi} :
A(\mathbb{D}) \rightarrow A$, via $z \mapsto \pi(e_1)$ and $ 1
\mapsto \pi(v_0)$.  By definition this yields a contractive
representation of $A(\mathbb{D})$ extending $\pi$.  Now as a
contractive representation of $A(\mathbb{D})$ is completely
contractive \cite[Corollary 3.14]{Muh-Sol:1998} the result is
established.
\end{example}
\begin{example}\label{arrow} Let $L_1$ be the finite graph given by
\[
\xymatrix{ {\bullet}_{v_0} \ar[r]^{e_1} & {\bullet}_{v_1} }
\]
We begin by letting $A_0$ be the nonunital subalgebra of
$A(\mathbb{D})$ generated by the coordinate function $z$. Consider
$A = \mathbb{C} * \mathbb{C} * A_0$ and let $p_0$ denote the image
of $v_0$ and $ p_1$ denote the image of $v_1$, under the
inclusions induced by sending one vertex to one copy of
$\mathbb{C}$ and the other vertex to the other copy of
$\mathbb{C}$. Letting $J$ be the ideal in $A$ generated by
elements of the form $zp_1 -z$ and $p_2z -z$, where $z$ is the
coordinate function in $\mathbb{C}$ then $OA(L_1) =
A/J$.\end{example}
\begin{remark} These examples serve to illustrate the difference between
the Toeplitz algebra of a quiver $T_+(Q)$ and the universal
operator algebra $OA(Q)$.  In fact, it will follow from our
analysis that $T_+(Q) = OA(Q)$ only for the graphs $V_0$ and
$B_1$.  Notice that $T_+(V_2) = \mathbb{C} \oplus \mathbb{C}$, and
$T_+(L_1) $ is the two by two upper triangular matrices.
\end{remark} These four examples will serve as building blocks for
all finite directed graphs.  We now look at how this is done.
First we will need some definitions and a lemma.
\begin{definition} We say that $k: Q \rightarrow R$ is a directed
graph homomorphism if $k = (k_E, k_V)$, an ordered pair of maps,
where $k_E: E(Q) \rightarrow E(R)$ and $k_V: V(E) \rightarrow
V(R)$ such that $s(k_E(e)) = k_V(s(e))$ and $ r(k_E(e)) =
k_V(r(e))$ for all $e \in E(Q)$.  We say that $k$ is a
monomorphism if $k_E$ and $k_V$ are both monomorphisms of
sets.\end{definition}
\begin{definition} We say that $Q$ is a directed subgraph of $R$, denoted $Q < R$,
if there exists a monomorphism $k: Q \rightarrow R$.
\end{definition}
We will often suppress mention of the monomorphism in what
follows.
\begin{lemma}\label{subalgebra} If $Q < R$, then $OA(Q) \subseteq OA(R)$.
\end{lemma}
\begin{proof} Let $ k: Q \rightarrow R$ be the monomorphism and
denote by $ \iota: R \rightarrow OA(R)$ the canonical contractive
representation given in Proposition \ref{oag}.  Notice that $\iota
\circ k: Q \rightarrow OA(R)$ induces a contractive representation
of $Q$. We let $S(Q)$ denote the subalgebra of $OA(R)$ generated
by $\iota \circ k (Q)$. We will show that $S(Q)$ satisfies the
requisite universal property of Proposition \ref{oag}.

Let $ \pi: Q \rightarrow A$ be a contractive representation of
$Q$.  Now define $p: R \rightarrow A$ via
\[ p(w) = \begin{cases} \pi(w) & \mbox{for } w \in k(Q) \\ 0 & \mbox{else}.
\end{cases}\]  Notice that $p$ is a contractive representation
of $R$.  There is then an extension $\tilde{p} : OA(R) \rightarrow
A$ such that $\tilde{p} \circ \iota(w) = p(w) $ for all $ w \in
w(R)$. It follows that $ \tilde{p}|_{S(Q)}$ is the requisite
extension and hence $S(Q)$ is $OA(Q)$. \end{proof}

We now define the free product of a graph and relate it to the
free product of operator algebras.
\begin{definition} Let $Q_1$ and $Q_2$ be two directed graphs.
Suppose $R$ is a directed subgraph of both $Q_1$ and $Q_2$; i.e.
$R < Q_1$ and $R < Q_2$. We define the vertex set of $Q_1*_RQ_2$
by \[ V(Q_1*_RQ_2) = (V(Q_1) \setminus V(R)) \sqcup (V(Q_2)
\setminus V(R)) \sqcup V(R). \]  We let \[ E(Q_1*_RQ_2) = (E(Q_1)
\setminus E(R)) \sqcup (E(Q_2) \setminus E(R)) \sqcup E(R)\] where
$\sqcup$ is defined to be disjoint union.  Notice that every $e
\in E(Q_1*_RQ_2)$ comes from a unique element of $E(Q_1)$ or
$E(Q_2)$ or perhaps both.  We define $r(e)$ to be the vertex in
$V(Q_1*_RQ_2)$ corresponding to the range in the original graph.
We define $s(e)$ similarly.
\end{definition}
\begin{remark} We have not required that $R$ be nonempty
and in fact we do not want to require this.  This allows us to
treat disconnected graphs as unamalgamated free products of the
connected components.\end{remark}
\begin{remark} Any finite graph can be constructed as
a finite free product of copies of $V_1$, $V_2$, $L_1$, and $B_1$.
This will make the next result central to our investigations.
\end{remark}

Recall the definition of the universal operator algebraic free
product of two operator algebras $A$ and $B$ \cite[Section
4]{Bl-Paul:1991} and its associated universal properties.  Care
must be exercised here, as $OA(Q)$ need not be unital. In fact we
must first adjoin a unit as in \cite{Mey:2001}. It is a corollary
of \cite[Theorem 3.2]{Mey:2001} that the free product of $A$ and
$B$ amalgamated over $C$, will be the subalgebra of
$A^+*_{C^+}B^+$ generated by $A$ and $B$.
\begin{theorem}\label{freeproducts} Let $Q_1$ and $Q_2$ be directed graphs with common
directed subgraph $R$.  Then $OA(Q_1*_RQ_2)$ is completely
isometrically isomorphic to the universal operator algebraic free
product $OA(Q_1)*_{OA(R)} OA(Q_2)$.\end{theorem}
\begin{proof} Notice by Lemma \ref{subalgebra}, that \[ OA(R) \subseteq
OA(Q_1) \mbox{ and } OA(R) \subseteq OA(Q_2) \] and hence
$OA(Q_1)*_{OA(R)} OA(Q_2)$ is well defined.  We will show that
$OA(Q_1)*_{OA(R)}OA(Q_2)$ satisfies the universal properties of
$OA(Q_1*_RQ_2)$.

Notice first, for $i = 1,2$ that there are contractive
representations\[ \pi_i: Q_i \rightarrow OA(Q_1)*_{OA(R)} OA(Q_2)
\] given by $\pi_i = \tau_i \circ \iota_i$ where $\tau_i: OA(Q_i)
\rightarrow OA(Q_1)*_{OA(R)} OA(Q_2)$ is the canonical inclusion,
and $\iota_i: Q_i \rightarrow OA(Q_i)$ is the canonical
contractive representation of Proposition \ref{oag}.  Notice that
$\pi_1|_R = \pi_2|_R$. It follows that there is a contractive
representation\[ \pi:Q_1*_RQ_2 \rightarrow OA(Q_1)*_{OA(R)}
OA(Q_2) \] and that $\pi(Q_1*_RQ_2) $ generates $OA(Q_1)*_{OA(R)}
OA(Q_2)$.

Now let $\sigma: Q_1*_RQ_2 \rightarrow A$ be a contractive
representation.  Then $\sigma|_{Q_1} =: \sigma_1$ and $
\sigma|_{Q_2} =: \sigma_2$ are contractive representations of
$Q_1$ and $Q_2$ respectively.  Furthermore, we know that
$\sigma_1(h) = \sigma_2(h)$ for all words $h \in w(R)$.  It
follows that there exists completely contractive homomorphisms
$\tilde{\sigma_1} : OA(Q_1) \rightarrow A$ and $\tilde{\sigma_2} :
OA(Q_2) \rightarrow A$, such that $\tilde{\sigma_1}|_{OA(R)} =
\widetilde{\sigma_2}|_{OA(R)}$. By the universal property for free
products of operator algebras there exists a completely
contractive homomorphism $\tilde{\sigma} : OA(Q_1)*_{OA(R)}OA(Q_2)
\rightarrow A$.  The result now follows.\end{proof}
\begin{example} We let $Q$ be the following graph:
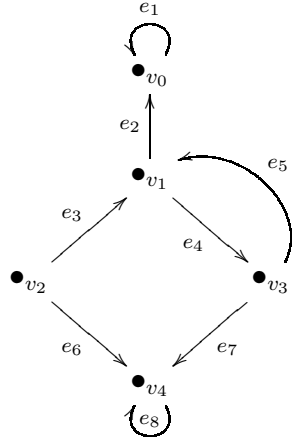
\[
\xymatrix{ & {\bullet}_{v_0} \ar@(ur,ul)[]_{e_1} & \\  &
{\bullet}_{v_1} \ar[u]^{e_2} \ar[dr]_{e_4}  & \\  {\bullet}_{v_2}
\ar[dr]_{e_6} \ar[ur]^{e_3} & & {\bullet}_{v_3} \ar[dl]^{e_7}
\ar@(ur,ur)[ul]_{e_5} \\ & {\bullet}_{v_4} \ar@(dr,dl)[]_{e_8}& }
\]
We decompose this graph into our building block graphs by first
looking at the fivefold unamalgamated free product of $V_0$ with
itself (i.e. the graph with five vertices and no edges).  We will
label the vertices with the labels $v_0, v_1, v_2, v_3,$ and
$v_4$.
\[
\xymatrix{ {\bullet}_{v_0} & {\bullet}_{v_1} & {\bullet}_{v_2} &
{\bullet}_{v_3}& {\bullet}_{v_4}}. \] We then amalgamate copies of
$L_1$ with the appropriate vertices (i.e. If \[ L_1 = \xymatrix{
{\bullet}_{u_0} \ar[r]^{t_1} & {\bullet}_{u_1} } \] then we
identify $u_0$ with $v_1$ and $u_1$ to $v_0$ to get the graph \[
\xymatrix{ {\bullet}_{v_0} & {\bullet}_{v_1} \ar[l]^{t_1} &
{\bullet}_{v_2} & {\bullet}_{v_3}& {\bullet}_{v_4}}). \] Lastly we
amalgamate copies of $B_1$ to $v_0$ and $v_4$. This yields the
decomposition
\[ \left[ (B_1 *_{v_0} L_1) *_{v_1} ((L_1 *_{v_2}L_1)*_{v_4}B_1) \right]
*_{\{ v_1, v_4 \}} \left[ L_1 *_{v_3}(L_1 *_{\{ v_1, v_3 \}}L_1)
\right] \]  and hence the operator algebra can be constructed
using free products.
\end{example}
\begin{remark}Notice that the decomposition of the graph is not unique, but
$OA(Q)$ is unique via universal considerations.\end{remark} Notice
that we must be careful of where the amalgamation is occurring.
For example the graphs
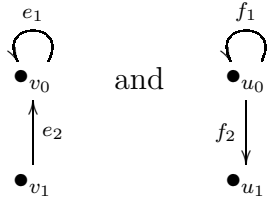
\[
\xymatrix{ {\bullet}_{v_0} \ar@(ur,ul)[]_{e_1} \\ {\bullet}_{v_1}
\ar[u]_{e_2}} \quad \mbox{ and } \quad \xymatrix{
{\bullet}_{u_0} \ar@(ur,ul)[]_{f_1} \ar[d]_{f_2} \\
{\bullet}_{u_1} }
\]
are similar but when constructing the free product the
amalgamation is being taken over different copies of $\mathbb{C}$
inside $OA(L_1)$.  In particular these two operator algebras are
adjoints of each other, and hence anti-isomorphic.  One might hope
for some sort of uniqueness as in \cite[Section 3]{Kat-Kribs:2003}
and \cite[Section 3]{Sol:2003}. We are currently investigating
such a result.

\section{The Universal $C^*$-Algebra of a Directed Graph}

We now want to look at a construction similar to that in Section
\ref{prelims}.  Here we will be concerned with building the
maximal $C^*$-algebra of a directed graph, as opposed to the
universal operator algebra of a directed graph. This is, in
spirit, very similar to a construction in Section 2 of
\cite{Good-Men:1990}. Before we actually construct the
$C^*$-algebra we calculate $OA(Q)^*$ for a directed graph $Q$.
\begin{definition} Let $Q$ be a directed graph.  We say that $Q^{\leftarrow}$
is the adjoint of $Q$ if $Q^{\leftarrow}$ is the graph obtained
from $Q$ by reversing all of the arrows. For notations sake, we
denote by $e^{\leftarrow}$ the reversed edge associated to the
edge $e$.\end{definition} The connection with adjoints is not
immediately obvious, but the next proposition justifies the name.
\begin{proposition} Let $Q$ be a directed graph, then $OA(Q)^* =
OA(Q^{\leftarrow})$. \end{proposition}
\begin{proof} Let $\pi: OA(Q) \rightarrow B(H)$ be a completely
isometric representation of $OA(Q)$.  Notice that this induces a
contractive representation $\pi: Q \rightarrow B(H)$. Now define a
contractive representation $\pi^*(Q^{\leftarrow})$ by $\pi^*(v) =
\pi(v)$ for all $ v \in V(Q)$ and $\pi^*(e^{\leftarrow}) =
\pi(e)^*$ for all $ e \in E(Q)$.  This induces a completely
contractive representation $\widetilde{\pi^*} : OA(Q^{\leftarrow})
\rightarrow B(H)$ such that $\widetilde{\pi^*}
(OA(Q^{\leftarrow})) = [\pi(OA(Q))]^*$.  The result now follows by
uniqueness of the algebra $[OA(Q)]^*$.\end{proof}

In section 2 we gave a method for constructing
\[ OA(Q)*_{OA(V(Q))}OA(Q)^*.\] In particular, we first construct the adjoint graph
and then look at the algebra $OA(Q *_{V(Q)} Q^{\leftarrow})$.  We
will use a similar construction to build the maximal $C^*$-algebra
of the graph.  The universal properties will be similar to those
of the maximal $C^*$-algebra of an operator algebra.  In fact we
will describe their relationship after the construction.
\begin{definition} Let $Q$ be a directed graph.  We say that $\pi:
Q*_{V(Q)}Q^{\leftarrow} \rightarrow C$ is a contractive $*$
representation of $Q$ if $\pi$ is a contractive representation
such that $ \pi(e)^* = \pi(e^{\leftarrow})$ for all $e \in
E(Q)$.\end{definition}
\begin{proposition}\label{cmg}For a directed graph $Q$ there exists a $C^*$-algebra,
$GC^*_m(Q)$ and a contractive $*$-representation $ \iota: Q
*_{V(Q)} Q^{\leftarrow} \rightarrow GC^*_m(Q)$ such that:
\begin{enumerate} \item $ \iota(Q *_{V(Q)} Q^{\leftarrow}) $ generates $GC^*_m(Q)$ as a
$C^*$-algebra. \item\label{cmgup} (Universal Property) for any
contractive $*$- representation $\varphi: Q *_{V(Q)}
Q^{\leftarrow} \rightarrow C$ into a $C^*$-algebra $C$ there is a
unique $*$-homomorphism $ \tilde{\varphi}: GC^*_m(Q) \rightarrow
C$ making the following diagram commute
\[
\xymatrix{ GC^*_m(Q) \ar@{-->}[dr]^{\tilde{\varphi}} & \\ Q
*_{V(Q)} Q^{\leftarrow} \ar[u]^{\iota} \ar[r]^{\varphi} & C }
\]
\end{enumerate}
\end{proposition}
\begin{remark} Notice that by the universal property, $GC^*_m(Q)$ is
unique up to $*$ isomorphism.\end{remark}

The construction follows a similar construction in Section 2 of
\cite{Good-Men:1990}.
\begin{proof} We proceed as in the construction of $OA(Q)$ but
this time we look at $\mathbb{C}w(Q*_{V(Q)}Q^{\leftarrow})$. We
define a seminorm on $\mathbb{C}w(Q*_{V(Q)}Q^{\leftarrow})$ by
setting \begin{align*} \| f \|_{C^*_m} = \sup \{ \| \pi(f) \| &:
\pi: Q *_{V(Q)} Q^{\leftarrow} \rightarrow C \mbox{ is a
contractive } \\ &* \mbox{-representation} \}. \end{align*} This
yields a $C^*$ seminorm on $\mathbb{C}w(Q*_{V(Q)}Q^{\leftarrow})$.
Denoting by $K$, the kernel of the seminorm, we have a norm on
$\mathbb{C}w(Q*_{V(Q)}Q^{\leftarrow}) / K$.  Completing with
respect to this norm yields a $C^*$-algebra satisfying the
universal properties.\end{proof}
\begin{remark}  Notice that $\| f \|_{C^*_m} \leq \| f
\|_{OA}$ as defined previously.  For those $f$ in the image of
$\iota: Q *_{V(Q)} Q^{\leftarrow} \rightarrow GC^*_m(Q)$, on the
other hand, $\| f \|_{C^*_m} = \| f \|_{OA}$.
\end{remark}
\begin{theorem} Let $Q$ be a directed graph.  Then $GC^*_m(Q)$ is
$C^*$ isomorphic to $C^*_m(OA(Q))$.\end{theorem}
\begin{proof} The proof is just repeated applications of universal
properties.  Notice that by the definition of $OA(Q)$ there exists
a completely contractive representation $\pi: OA(Q) \rightarrow
GC^*_m(Q)$ which extends the map $\iota|_Q$.  By the definition of
$C^*_m(OA(Q))$ it follows that there is a $*$ homomorphism
$\tilde{\pi}: C^*_m(OA(Q)) \rightarrow GC^*_m(Q)$.

There is a canonical contractive representation $\sigma: Q
\rightarrow C^*_m(OA(Q))$.  Extending this to a $*$ representation
$\tilde{\sigma}: GC^*_m(Q) \rightarrow C^*_m(OA(Q))$ we need only
verify that the two $*$ representations are inverses of each
other.  Notice that $ \widetilde{\sigma}(w) = \pi(w)$ for $w \in
V(Q) \cup E(Q)$ by definition.  Similarly $\widetilde{\pi}(w) =
\sigma(w)$ for $ w \in V(Q) \cup E(Q)$.  The result now
follows.\end{proof}

Blecher \cite[Proposition 2.2]{Bl:1999} observed that for operator
algebras $A$ and $B$, $C^*_m(A*B) = C^*_m(A) * C^*_m(B)$ where
$A*B$ is universal in the category of operator algebras with
completely contractive homomorphisms and $C^*_m(A) * C^*_m(B)$ is
universal in the category of $C^*$-algebras with $*$
homomorphisms. This result allows us to build $GC^*_m(Q)$ for
finite graphs $Q$ in a manner similar to the way we built $OA(Q)$.
\begin{proposition} Let $Q_1$ and $Q_2$ be finite graphs
and let $R$ be some (sub)collection of vertices in $V(Q_1)\cap
V(Q_2)$. Then $GC^*_m(Q_1 *_RQ_2)$ is $C^*$ isomorphic to
\[ GC^*_m(Q_1) *_{\mathbb{C} * \mathbb{C} \cdots *
\mathbb{C}} GC^*_m(Q_2) \] where there are $|R|$ copies of
$\mathbb{C}$ in the preceding formula.\end{proposition}

We remind the reader that we do not assume that a contractive $*$
representation of a directed graph send distinct vertices to
orthogonal projections. Nor for that matter do we need assume that
edges be sent to partial isometries.  Notice that even in the
simple case of graphs with multiple vertices and no edges we do
not have the universal $C^*$-algebra $C^*(Q)$ of Pask, Raeburn,
Renault, etc. On the other hand, notice that the universal
$C^*$-algebra will yield a contractive $*$ representation of $Q$.
It follows from the universal properties of $GC^*_m(Q)$ that
$C^*(Q)$ is a quotient of $GC^*_m(Q)$.

\section{Continuity of universal algebras with respect to direct limits}

So far, we have not needed to restrict to finite graphs, although
our examples have been using finite directed graphs. In this
section we will give a method for dealing with infinite graphs. We
begin by looking at inductive limits of operator algebras. This
approach was inspired by a result of Spielberg \cite[Theorem
2.35]{Spielberg:2002}.

We say that $(A_i, \varphi_{ij})$ is an inductive system of
operator algebras if each $A_i$ is an operator algebra and each $
\varphi_{ij} : A_j \rightarrow A_i$ is a completely contractive
homomorphism such that $ \varphi_{ij} \circ \varphi_{jk} =
\varphi_{ik}$.  Notice that if each $A_i$ is a $C^*$-algebra then,
since a completely contractive homomorphism is a $*$ homomorphism,
this is the usual definition of an inductive system of
$C^*$-algebras. Recall, as in \cite[Appendix L]{Wegge:1993}, that
to an inductive system $(A_i, \varphi_{ij})$ of operator algebras
there exists a unique inductive limit operator algebra $A=
\displaystyle{ \lim_{\rightarrow} A_i}$ satisfying the following
universal property:
\begin{quote}\label{ilup} (Universal Property) There exist
completely contractive homomorphisms $ \Phi_i : A_i \rightarrow A$
such that $ \Phi_i \circ \varphi_{ij}(x)  = \Phi_j (x)$ for all $x
\in A_j$.  Further if $B$ is an operator algebra such that there
are completely contractive homomorphisms $\Gamma_i: A_i
\rightarrow B$ with $ \Gamma_i \circ \varphi_{ij}(x)  = \Gamma_j
(x)$ for all $x \in A_j$, then there is a unique completely
contractive homomorphism $ \Lambda : A \rightarrow B$ such that $
\Lambda \circ \Phi_i (x) = \Gamma_i(x)$ for all $x \in A_i$ so
that the following diagram commutes:
\[
\xymatrix{
A_j \ar[r]^{\Phi_j} \ar[dr]|<<{\Gamma_j} & A \ar@{-->}[d]^{\Lambda} \\
A_i \ar[u]^{\varphi_{ji}} \ar[ur]|<<{\Phi_i} \ar[r]_{\Gamma_i} &
B}
\]
\end{quote}
\begin{theorem} \label{cmil} Let $(A_i, \varphi_{ij})$ be
an inductive system of operator algebras.  Then there exist maps $
\tilde{\varphi}_{ij}: C^*_m(A_j) \rightarrow C^*_m(A_i)$ such that
$ (C^*_m(A_i), \tilde{\varphi}_{ij})$ is an inductive system of
$C^*$-algebras.  Furthermore \[ \lim_{\rightarrow} C^*_m(A_i) =
C^*_m( \lim_{\rightarrow}A_i).\]
\end{theorem}
\begin{proof}  By the universal property of $C^*_m(A_i)$, each $
\varphi_{ij}: A_j \rightarrow A_i \subseteq C^*_m(A_i)$ will
extend to a $*$ homomorphism $ \tilde{\varphi}_{ij}: C^*_m(A_j)
\rightarrow C^*_m(A_i) $. By uniqueness of such extensions we get
that $ \tilde{\varphi}_{ij} \circ \tilde{\varphi}_{jk} =
\tilde{\varphi}_{ik}$.  Thus $(C^*_m(A_i), \tilde{\varphi}_{ij})$
is an inductive system.

Now let $A = \displaystyle{ \lim_{\rightarrow}} A_i$ and $C =
\displaystyle{ \lim_{\rightarrow}} C^*_m(A_i)$.  We need to show
that $C^*_m(A)$ satisfies the universal property unique to $C$.
Assume that $B$ is a $C^*$-algebra, and $ \Gamma_{i} : C^*_m(A_i)
\rightarrow B$ are $*$ homomorphisms.  Further assume that
$\sigma_i \circ \tilde{\varphi}_{ij}(x)  = \sigma_j (x) $ for all
$ x \in C^*_m(A_j)$.  Let $ \overline{\sigma_j} = \sigma_j|_{A_j}
$. Notice that $\overline{\sigma_j}$ is completely contractive and
$ \overline{\sigma_i} \circ \varphi_{ij}(x)  =
\overline{\sigma_j}(x) $ for all $ x \in A_j$. Now by the
universal property for $A$, there is a completely contractive
homomorphism $ \overline{\sigma} : A \rightarrow B $ with $
\overline{\sigma} \circ \Phi_j (x) = \overline{\sigma_j} (x) $ for
all $ x \in A_j$.  The universal property for $C^*_m(A)$ now gives
a $*$ homomorphism $ \sigma: C^*_m(A) \rightarrow B$ satisfying $
\sigma \circ \tilde{\Phi}_j (x) = \sigma_j (x) $ for $x \in
C^*_m(A_j)$.  By uniqueness of the inductive limit it follows that
$\displaystyle{C^*_m(A) \cong \lim_{\rightarrow} C^*_m(A_i)}$.
\end{proof}

Now we look at a similar result in the context of directed graphs.
\begin{theorem}\label{oagil} Let $Q$ be a directed graph and let $W$ be any
collection of subgraphs such that
\[ \cup_{F \in W} V(F) = V(Q) \] and \[ \cup_{F \in W} E(F) =
E(Q). \] Then \[ OA(Q) = \lim_{\rightarrow} OA(W).\]\end{theorem}
\begin{proof} We show that $OA(Q)$ has the requisite universal
property.  First notice that by Lemma \ref{subalgebra} it follows
that there exists completely contractive maps $\Gamma_R: OA(R)
\rightarrow OA(Q)$ for all $R \in F(Q)$.  Further if $ R \subseteq
S \subseteq Q$ it follows that there is a connecting completely
contractive homomorphism $\Gamma_{SR}: OA(R) \rightarrow OA(S)$.
Now for each $R \in F(Q)$, let $\pi_R: OA(R) \rightarrow A$ be a
completely contractive representation such that $\pi_S \circ
\Gamma_{SR}|_{OA(R)} = \pi_R|_{OA(R)}$. Then define $\pi: Q
\rightarrow A$ via $\pi(e) = \pi_{E}$, where $E$ is the subgraph
given by $\{ e, r(e), s(e) \}$. Now by definition $\Gamma_R(OA(R))
= \Gamma_S \circ \Gamma_{SR} (OA(R))$ hence it follows that
$OA(Q)$ is in fact the inductive limit.\end{proof}
\begin{corollary}\label{cmgil} Let $Q$ be a directed graph and $W$
as in the previous theorem.  Then $\displaystyle{GC^*_m(Q) =
\lim_{\rightarrow}C^*_m(Q_i)}$. \end{corollary}
\begin{proof}This follows directly using Theorems \ref{cmil} and \ref{oagil}.\end{proof}

We can now construct operator algebras for a graph $Q$ with
countable edge and vertex sets.  We first construct the operator
algebras of the finite directed subgraphs as free products of
$V_0$, $L_1$ and $B_1$. We then use inductive limits to concretely
build both $OA(Q)$ and $GC^*_m(Q)$. We will apply this technique
in the next section when we discuss the K-Theory of universal
operator algebras of a directed graph.
\begin{remark} Notice that as $OA(V_2)$, $OA(L_1)$, and
$A(\mathbb{D})$ are not finite dimensional it follows that $OA(Q)$
is not the direct limit of finite dimensional operator algebras.
On the other hand, the result of Spielberg \cite[Theorem
2.35]{Spielberg:2002} suggest that for the Toeplitz algebra of a
quiver of Muhly \cite[Section 6]{Muhly:1997} continuity with
respect to direct limits may yield $AF$ algebras in certain
cases.\end{remark}

\section{K-Theory}

We now use a generalized version of a result of Cuntz
\cite{Cu:1982} to allow us to calculate the K-Theory using the
free product decomposition of a finite directed graph.
\begin{theorem}\label{cuntz}Let $A$ and $B$ be operator algebras sharing a
common $C^*$ subalgebra $D$. Further assume that there exists onto
idempotent completely contractive homomorphisms $ \varphi_A : A
\rightarrow D$ and $ \varphi_B: B \rightarrow D$. Let $ j_1: D
\rightarrow A$, $j_2: D \rightarrow B$, $i_1 : A \rightarrow
A*_DB$, and $ i_2 : A \rightarrow A*_DB $ be the canonical
inclusions.

There is a natural split exact sequence of K-groups \[ \xymatrix{
0 \ar[r] & K_*(D) \ar[r]^-{d} & K_*(A) \oplus K_*(B) \ar[r]^-{i} &
K_*(A*_DB) \ar[r] & 0 } \] where \[ d(x) = (j_{1*}(x), -j_{2*}(x))
\quad \mbox{ and } \quad i(x) = i_{1*}(x) + i_{2*}(x).
\]
\end{theorem}
\begin{proof} The proof of Theorem \ref{cuntz} follows in the same manner as the
corresponding result in \cite{Cu:1982}.  We do not give the
details here.\end{proof} This allows us to compute the K-groups
for free products in a rather restricted case. Fortunately, our
universal graph algebras are a class of operator algebras with
idempotent homomorphisms onto common subalgebras. In effect, then,
this result provides a method by which we can construct the
K-groups for an arbitrary directed graph with countable edge and
vertex set. The actual group will have to be done on a case by
case basis, but the general idea is as follows:

We first look at finite subgraphs of $Q$.  For a finite subgraph
$Q'$ we can use the free product representation to decompose
$OA(Q')$ as a finite number of free products of $OA(V_1)$,
$OA(V_2)$, $OA(L_1)$, and $OA(B_1)$. Noting that the projection
onto the vertices induces a completely contractive idempotent
representation of $OA(Q')$ onto $OA(V(Q'))$ we can apply the above
result to calculate $K_*(OA(Q'))$.  Then using inductive limits we
can actually calculate $\displaystyle{ K_*(OA(Q)) =
\lim_{\rightarrow} K_*(OA(F(Q)))}$, where $F(Q)$ is the set of
finite subgraphs ordered by inclusion.

The following result allows us actual computation of the K-groups
for finite graphs $Q$.  We will also use the following result to
calculate the K-groups for $GC^*_m(Q)$.
\begin{proposition}\label{kfingraph} Let $Q$ be a finite directed
graph. Then there exists a family of completely contractive
homomorphisms $ \varphi_t: OA(Q) \rightarrow OA(Q)$ such that
$\varphi_1 = id_{OA(Q)}$ and the range of $\varphi_0$ is
$OA(V(Q))$.  Moreover, $\{ \varphi_t \}_{t \in [0,1]}$ is
pointwise norm continuous.\end{proposition}
\begin{proof} We begin by defining $ \varphi_t$ on a dense
subalgebra of $OA(Q)$.  Letting $ \iota: Q \rightarrow OA(Q)$ be
the canonical inclusion, define the map $ \varphi_t: Q \rightarrow
OA(Q)$ given by $\varphi_t(e) = t \iota(e)$ for $e \in E(Q)$ and
$\varphi_t(v) = \iota(v)$ for $v \in V(Q)$.  Clearly $\varphi_t$
is a contractive representation for all $ 0 \leq t \leq 1$. Thus
by the universal property there exists a completely contractive
homomorphism $ \tilde{\varphi_t} : OA(Q) \rightarrow OA(Q)$
extending $ \varphi_t$.  Clearly $\tilde{\varphi_1} = id_{OA(Q)}$
and the range of $ \tilde{\varphi_0}$ is $OA(V(Q)) \subseteq
OA(Q)$.

We need only show that $ \tilde{\varphi_t}$ is pointwise norm
continuous. Hence let $f \in OA(Q)$, and let $ \varepsilon > 0$ be
given.  First notice that as $\mathbb{C}(w(Q))$ is dense in
$OA(Q)$ there exists some element $f_1 \neq 0 \in \mathbb{C}w(Q)$
such that $\| f_1 - f \| < \frac{\varepsilon}{3}$.

Let $n$ be the length of the longest word in $w(Q)$ such that $
f_1(w) \neq 0$.  Let $m$ be the cardinality of the set of $W = \{
w \in w(Q)\setminus V(Q): f_1(w) \neq 0\}$.  Now let $M = n  m \|
f_1 |_{W} \| $ and $ \delta = \frac{\varepsilon}{3(M+1)}$. We
claim that if $ |t-t_0|< \delta$ then $ \| \varphi_t(f) -
\varphi_{t_0} (f) \| < \varepsilon $.

Remembering that $|w|$ is the length of a word in $w(Q)$, we
notice:
\begin{align*} \| \varphi_t(f) -& \varphi_{t_0} (f ) \| \\ & \leq
\| \varphi_t(f) -
\varphi_t(f_1) \| + \| \varphi_{t}(f_1) - \varphi_{t_0}(f_1) \| \\
 & \hspace{.5cm} + \| \varphi_{t_0} (f_1) - \varphi_{t_0}(f) \| \\ & \leq  \| f -
f_1 \| + \| \varphi_t(f_1) - \varphi_{t_0}(f_1) \| + \| f- f_1 \|
\\ & <  \frac{2 \varepsilon}{3} + \| \varphi_t(f_1) -
\varphi_{t_0}(f_1) \| \\ & \leq  \frac{2 \varepsilon}{3} +  \|
\sum_{w \in W} (t^{|w|} - t_0^{|w|}) f(w) \| \\ & <  \frac{2
\varepsilon}{3} \\ & \hspace{.5cm}+ | t- t_0|\cdot
\| \sum_{w \in W} (t^{|w|-1} + t^{|w|-2}t_0 + \cdots + t_0^{|w|-1})f(w) \| \\
& <  \frac{2 \varepsilon}{3} + \delta \| \sum_{w \in W} n f(w) \|
\\ & < \varepsilon. \end{align*} The result now follows.
\end{proof}

Hence by the homotopy invariance of $K$-groups \[K_*(OA(Q)) =
K_*(OA(V(Q))).\] Thus an application of Theorem \ref{cuntz},
yields
\[ K_*(OA(Q)) = K_*(\mathbb{C} \oplus \mathbb{C} \oplus \cdots
\oplus \mathbb{C})
\] where there are $|V(Q)|$ copies of $\mathbb{C}$ in the previous
formula.
\begin{remark} Notice that the K-groups do not constitute a very
good invariant since, in effect, it merely keeps track of the
vertices, while all edge data is lost.\end{remark}
\begin{remark}  The above calculations can be adapted for the
Toeplitz algebra of a quiver and hence it can be shown that
$K_*(T_+(Q)) = K_*(OA(Q))$. \end{remark}

We now look at some results for $C^*_m(A)$ where $A$ is an
operator algebra. These results will allow us to compute the
K-groups for $GC^*_m(Q)$.  To accomplish this we discuss how
completely contractive endomorphisms of $A$ extend to
$*$-endomorphisms of $C^*_m(A)$ and some properties preserved by
the extension.  A completely contractive endomorphism $\varphi: A
\rightarrow A$ can be treated as a map $\varphi: A \rightarrow
C^*_m(A)$ and hence there is an extension $ \tilde{\varphi} :
C^*_m(A) \rightarrow C^*_m(A)$ which is a $*$ endomorphism of
$C^*_m(A)$. We now start with a simple observation.
\begin{proposition} \label{homotopy} Let $\{ \varphi_t \}_{t \in [0,1]} $ be a set
of completely contractive endomorphisms of A which is point-norm
continuous with respect to $t$.  Then $\{ \tilde{ \varphi}_t \}_{t
\in [0,1]} $ is a set of endomorphisms of $C^*_m(A)$ which is also
point-norm continuous with respect to $t$. \end{proposition}
\begin{proof}  We will show that for a fixed $a \in C^*_m(A)$ we
have that $ \tilde{\varphi}_t(a)$ is norm continuous with respect
to $t$. Recall that $P = A *_{a}A^* \subseteq C^*_m(A)$ is a dense
set in $C^*_m(A)$. Notice also that continuity of $ {\varphi_t}^*$
is clear.

We first deal with elements of $P$ of the form
\begin{equation} \label{products} a = a_1 * a_2 * \cdots * a_m.
\end{equation}
We proceed by induction on the length of the word $a$.  If
$\ell(a) =1 $ then by hypothesis $ \varphi_t(a)$ is norm
continuous with respect to $t$.  Assume that $ \tilde{\varphi}_t
(a)$ is norm continuous for words of length $m-1$.  Letting $a =
a_1 * a_2 * \cdots * a_m$ then we can rewrite $a$ as $a_1 * b$
where $ b = a_2 * a_3 \cdots *a_m$.  By the induction hypothesis,
given $\varepsilon > 0 $ there exists $ \delta > 0$ such that for
$ |t - t_0| < \delta $ we have $\| \tilde{\varphi}_t(a_1) -
\tilde{\varphi}_{t_0} (a_1) \| < \dfrac{\varepsilon}{2 (\|
\tilde{\varphi}_{t_0}(b) \| + 1)} $ and $\| \tilde{\varphi}_t(b) -
\tilde{\varphi}_{t_0} (b) \| < \dfrac{\varepsilon}{2 (\|
\tilde{\varphi}_{t}(a_1) \| + 1)}. $ Now for $| t- t_0| < \delta$
\begin{eqnarray*} \| \tilde{\varphi}_t(a_1 * b) -
\tilde{\varphi}_{t_0} (a_1 * b)  \| & = & \| \tilde{\varphi}_t(a_1
*b) - \tilde{\varphi}_{t}(a_1)\tilde{\varphi}_{t_0}(b) \\ & &
\quad + \tilde{\varphi}_{t}(a_1)\tilde{\varphi}_{t_0}(b) -
\tilde{\varphi}_{t_0}(a_1* b) \| \\ & \leq & \|
\tilde{\varphi}_t(a_1 *b) -
\tilde{\varphi}_{t}(a_1)\tilde{\varphi}_{t_0}(b) \| \\ & & \quad +
\| \tilde{\varphi}_{t}(a_1)\tilde{\varphi}_{t_0}(b) -
\tilde{\varphi}_{t_0}(a_1* b) \| \\ & = & \|
\tilde{\varphi}_{t}(a_1)\| \| \tilde{\varphi}_{t_0}(b) -
\tilde{\varphi}_t(b) \| \\ & & \quad + \|\tilde{\varphi}_t( a_1 )
- \tilde{\varphi}_{t_0}(a_1) \| \|\tilde{\varphi}_{t_0}(b) \| \\ &
< & \varepsilon. \end{eqnarray*}

We now establish continuity for all of $P$.  Linearity implies
that $ \tilde{\varphi}_t$ is norm continuous for $x \in P$.  We
can now show continuity for all of $C^*_m(A)$ since $P$ is dense
in $C^*_m(A)$. Let $ b$ be an arbitrary element of $C^*_m(A)$. For
$ \varepsilon
> 0 $ choose $ b_0 \in P$ such that $ \| b - b_0 \| <
\frac{\varepsilon}{3}$, and notice that $ \|
\tilde{\varphi_t}(b-b_0) \| < \frac{\varepsilon}{3} $  for all $t$
as $ \tilde{\varphi_t}$ is completely contractive for all $t$. Now
choose $ \delta > 0 $ as above so that $ | t- t_0 | < \delta $
implies $ \| \tilde{\varphi_t} (b_0) - \tilde{ \varphi_{t_0}} (
b_0) \|$.  Now notice that
\begin{eqnarray*} \| \tilde{\varphi_t}
( b ) - \tilde{\varphi_{t_0}} ( b ) \| & \leq & \|
\tilde{\varphi_t} ( b ) - \tilde{\varphi_{t}} ( b_0 ) \| + \|
\tilde{\varphi_t} ( b_0 ) - \tilde{\varphi_{t_0}} ( b_0 ) \| \\ &
& \quad  + \|
\tilde{\varphi_{t_0}} ( b_0 ) - \tilde{\varphi_{t_0}} ( b ) \| \\
& < & \frac{\varepsilon}{3} + \frac{ \varepsilon}{3} + \frac{
\varepsilon}{3} \\ & = & \varepsilon. \end{eqnarray*} \end{proof}

As an application it follows that $K_*(OA(B_1)) = K_*(C^*_m(B_1))
= K_*(\mathbb{C})$ and also $K_*(OA(L_1)) = K_*(C^*_m(L_1)) =
K_*(\mathbb{C} \oplus \mathbb{C})$. We then have the following
corollaries.
\begin{corollary}Let $Q$ be a finite directed graph, then $K_0(OA(Q)) =
K_0(GC^*_m(Q))$ and $ K_1(OA(Q)) = K_1(GC^*_m(Q))$.\end{corollary}
\begin{proof}First decompose $Q$ as a free product of copies of
$L_1$ and $B_1$.  As $K_*(OA(B_1)) = K_*(C^*_m(B_1))$ and
$K_*(OA(L_1)) = K_*(C^*_m(L_1))$, we can use Theorem
\ref{kfingraph} to show that $K_*(GC^*_m(Q)) = K_*(OA(Q))$.
\end{proof}
\begin{corollary} Let $Q$ be a directed graph with countable edge
and vertex sets. Then $K_*(OA(Q)) =
K_*(GC^*_m(Q))$.\end{corollary}
\begin{proof} This follows after applying theorem \ref{oagil}
and corollary \ref{cmgil}, as K-groups are continuous with respect
to direct limits.
\end{proof}

This suggests the following question:
\begin{question} For which operator algebras $A$ is it true that
$K_*(C^*_m(A)) = K_*(A)$? \end{question} From proposition
\ref{homotopy}, any operator algebra which has a pointwise norm
continuous homotopy onto a $C^*$ subalgebra will provide an
example where the groups coincide.  Hence for a counterexample one
would need a non-self adjoint operator algebra which does not have
a pointwise norm continuous homotopy onto the diagonal.

\begin{acknowledgement}The author would like to thank the referee
for many helpful comments.  This work comprises part of the
author's doctoral dissertation and he gratefully acknowledges the
support of his department and advisor, David
Pitts.\end{acknowledgement}

\bibliographystyle{amsplain}
\bibliography{common}

\end{document}